\documentclass[11pt]{article}
\usepackage{amssymb,amsmath}

\def\={~=~}
\newcommand{\lp}{\left(}
\newcommand{\lb}{\left\lbrack}
\newcommand{\rp}{\right)}
\newcommand{\rb}{\right\rbrack}

\newcommand{\Z}{\mathbb{Z}}

\begin{document}

\begin{center}
\Large
On the Use of Integrals to Evaluate  Series of Rational Terms
\end{center}

\vspace*{0.2 in}

\begin{flushright}
Costas J. Efthimiou  \\
Department of Physics \\
University of Central Florida \\
Orlando, FL 32816  \\
USA
\end{flushright}

\paragraph*{Introduction}
In the January 2006 issue of \textit{The College Mathematics
Journal}, M. Andreoli posed the following problem \cite{Andreoli}:

\fbox{\vbox{
\begin{quote}
 \textbf{Problem 819}. For $n\ge 1$, evaluate
 \begin{equation}
 \sum_{k=1}^{\infty}  {1\over k(k+1)(k+2)\cdots(k+n)}~.
 \label{eq:problem}
 \end{equation}
\end{quote}
    }}

This article is written with the hope to draw attention to a
method that has been proposed independently by the author in
\cite{Efthimiou} and Wheelon in \cite{Wheelon} and that allows one
to find exact values for a large class of convergent series of
rational terms. The method has been extended to many series in the
papers of Lesko and Smith \cite{Lesko} and Efthimiou
\cite{Efthimiou2}.

We will outline both the author's and Wheelon's variations of the
method in order to compare the similarities and differences, and
then use them to compute the series (\ref{eq:problem}). Finally,
we apply the methods to additional series that generalize
(\ref{eq:problem}). The results thus obtained seem to be new.

\paragraph*{The Core Idea of the Method}
Let a  series of the form $\sum_{k\in I} u_k v_k$, where $I$ is a
subset of $\Z$, be given. Then it is  convenient to write only one
of the factors, say $v_k$, as an integral transform
$$
    v_k \= \int_a^b \tilde v_k(t) \, f(t)\, dt~,
$$
and thus
\begin{eqnarray*}
  \sum_{k\in I} u_k v_k
         &=& \sum_{k\in I}u_k\,\int_a^b \tilde v_k(t) \,f(t)\,dt~.
\end{eqnarray*}
Assuming that the order of the operations of summation and integration
can be exchanged
\begin{eqnarray*}
  \sum_{k\in I} u_k v_k
         &=& \int_a^b\lp\sum_{k\in I}u_k \,\tilde v_k(t)\rp\,f(t)\,dt~.
\end{eqnarray*}
In this article we shall always exchange the order of the two
operations assuming that the reader knows how to reason for its
validity. Details on this may be found in the original papers
\cite{Efthimiou} and \cite{Lesko} where the Laplace transform is
used. If one can find an explicit function $h(t)=\sum_{k\in
I}u_k\,\tilde v_k(t)$, then he has succeeded in writing the
initial series in a simple integral representation:
\begin{eqnarray*}
    \sum_{k\in I} u_k v_k &=& \int_a^b h(t)\,f(t)\,dt~.
\end{eqnarray*}
If, furthermore, the integration can be performed, then analytic
answers for the initial series are obtained.

\paragraph*{Series of rational terms}
Consider a series for which $I=\{1,2,\dots\}$, $u_k=1$, and
$v_k=Q(k)/P(k)$, where $Q(k)$ and $P(k)$ are two polynomials in
$k$:
\begin{equation}
  S\=\sum_{k=1}^{\infty}\,{Q(k)\over P(k)}~.
 \label{eq:2}
 \end{equation}
 To ensure convergence, we may assume that deg$Q+2\le$deg$P$.

Problem 819 of \emph{The College Mathematics Journal} is clearly a
special case of (\ref{eq:2}) for $Q(k)=1$ and $P(k) \=
k(k+1)(k+2)\cdots(k+n)$.

\paragraph*{The author's variation}
 In the method presented in \cite{Efthimiou},
the expansion of the general rational term of the series to
partial fractions plays a central role. We shall then assume that
the sum $S$ can be written in terms of partial fractions:
$$
  S \=\sum_{k=1}^\infty\sum_{i=1}^\ell\sum_{j=1}^{m_i}\,
           {A_{ij}\over (k+a_i)^j}~,
$$
where the constants $A_{ij}$ are uniquely determined by the
partial fraction decomposition of each summand. Details for this
decomposition can be found in \cite{Efthimiou}.

Using the identity
\begin{equation}
\label{laplace2}
   {1\over A^L}\={1\over (L-1)!}\,
   \int_0^{+\infty}\, x^{L-1} e^{-Ax}\, dx~,
\end{equation}
we write the series in integral form (which is valid only if
$a_i>-1,$ for all $i$):
\begin{eqnarray}
  S &=&\sum_{k=1}^{\infty} \sum_{i=1}^\ell \sum_{j=1}^{m_i}\,
  {A_{ij}\over (j-1)!}\, \int_0^{+\infty}\,
    x^{j-1}\, e^{-(k+a_i)x}\, dx ~.
\end{eqnarray}
By interchanging the integration over $x$ with the addition over
$k$, we find
\begin{eqnarray}
  S &=&\sum_{i=1}^\ell \sum_{j=2}^{m_i}\,
  {A_{ij}\over (j-1)!}\,
   \int_0^{+\infty}\,
    x^{j-1}\, {e^{-(a_i+1)x}\over 1- e^{-x} }\,
    \, dx
  \nonumber \\
  &+&  \sum_{i=1}^\ell \,
  A_{i1}\,
   \int_0^{+\infty}\,
    {e^{-(a_i+1)x}-1\over 1- e^{-x} }\,
    \, dx~.
   \nonumber
\end{eqnarray}
To further simplify this sum we need to know the values of the
constants $A_{ij}$. The integrals may be expressed in terms of the
polygamma functions \cite{Efthimiou}.

\paragraph*{Andreoli's problem}     In problem 819 of \emph{The College
Mathematics Journal} the rational function $1/P(k)$ can be easily
expanded in partial fractions
 \begin{equation}
  {1\over P(k)}\= \sum_{i=0}^n {A_i\over k+i}~,
 \label{eq:1}
 \end{equation}
with
$$
   A_i \= {(-1)^i\over n!} \, \lp\begin{matrix}{n}\\{i}\end{matrix}\rp ~.
$$
This can be proved easily by multiplying equation (\ref{eq:1}) by
$k+j$ and then taking the limit $k\to-j$. We thus have
 \begin{eqnarray}
  \sum_{k=1}^\infty {1\over P(k)}
  &=& {1\over n!} \, \int_0^{+\infty} {e^{-x}\over 1-e^{-x}}
      \sum_{i=0}^n \lp\begin{matrix}{n}\\{i}\end{matrix}\rp
      (-1)^{ix} (e^{-x})^i \, dx \nonumber \\
  &=& {1\over n!}\, \int_0^{+\infty}  {e^{-x}\over 1-e^{-x}}
  (1-e^{-x})^n \, dx \nonumber \\
  &=& {1\over n!}\, \int_0^{+\infty}  e^{-x}
  (1-e^{-x})^{n-1} \, dx ~.
  \label{eq:10}
\end{eqnarray}
By a change of variables $u=e^{-x}$,
 \begin{eqnarray}
  \sum_{k=1}^\infty {1\over P(k)} \= {1\over n!}\int_0^1
  (1-u)^{n-1} du\= {1\over n \, n!}~.
 \end{eqnarray}

\paragraph*{A generalization}
 One can easily generalize Andreoli's proposal:
  For $n\ge 1$, we can evaluate
 \begin{equation}
 S(a,b)\=\sum_{k=1}^{\infty}  {1\over [a+kb][a+(k+1)b][a+(k+2)b]\cdots[a+(k+n)b]}~.
 \label{eq:9}
 \end{equation}
Following the method described above, we find
 $$
  S(a,b) \= {1\over b^n n!} \, \int_0^{+\infty}
  e^{-(a+b)x}(1-e^{-bx})^{n-1} \, dx~.
 $$
 The change of variables $u=e^{-x}$ gives the integral the form
 $$
  S(a,b) \= {1\over b^n n!} \, \int_0^1
  u^{a+b-1}(1-u^b)^{n-1} \, du~.
 $$
 The last integral is a  B-function integral:
 \begin{equation}
  \int_0^1  \, (1-u^\ell)^{n-1} \, u^{m-1}\, du \= {1\over\ell}\,
  B({m\over\ell},n)~,~~~~~m,n,\ell>0~.
\label{eq:beta}
\end{equation}
 The sum $S(a,b)$ is thus equal to $B(a/b+1,n)$ which value can be easily computed through
the $\Gamma$-function. Finally, we get
 $$
   S(a,b) \= {1\over n b} \, \prod_{i=1}^n {1\over a+i b}~.
 $$

\paragraph*{Wheelon's variation}

\def\buildrel#1\below#2{{\mathrel{\mathop{\kern0pt #2}\limits_{{\footnotesize #1}}}}}

The method presented by Wheelon  \cite{Wheelon} for the
computation of the series (\ref{eq:2}) assumes that $Q(k)=1$ and
$P(k)$ can be factorized as a product of of first order monomials.
In this method,  the convolution theorem as used by Feynman
\cite{Feynman} plays a central role. In particular, starting from
$$
  {1\over ab} \= \int_0^1 {dx \over [ax+b(1-x)]^2}~,
$$
one can show the general formula
\begin{eqnarray*}
 {1\over \prod\limits_{i=1}^\ell \alpha_i^{m_i}  } \=
 {\Gamma(\sum\limits_{i=1}^\ell m_i)\over
  \prod\limits_{i=1}^\ell \Gamma(m_i)}
     {\buildrel \begin{matrix}{0\le \sum\limits_{i=1}^{\ell-1}x_i \le1}\\
                           {0\le x_j \le 1}
              \end{matrix}
   \below {\int\int\cdots\int} }
   {  (\prod\limits_{i=1}^{\ell-1} x_i^{m_i-i}) \,
    (1-\sum\limits_{i=1}^{\ell-1} x_i)^{m_\ell-1} \over
    [\sum\limits_{i=1}^{\ell-1} x_i\alpha_i +
     (1-\sum\limits_{i=1}^{\ell-1} x_i)\alpha_\ell]^{\sum\limits_{i=1}^\ell m_i} }
     \, dx_1 dx_2 \cdots dx_{\ell-1} ~.
\end{eqnarray*}
In physics literature, this formula is known as Feynman's integral
and the $x$'s are called Feynman's parameters.

If $Q(k)=1$ and $P(k)$ given by
\begin{equation}
  P(k)= (k+a_1)^{m_1}(k+a_2)^{m_2}\dots (k+a_\ell)^{m_\ell}~,
 \label{eq:3}
 \end{equation}
where the $a_i,~i=1,2,\dots,\ell$ are distinct real numbers, none
of them a negative integer and all $m_i$ are positive integers,
then the series (\ref{eq:2}) is written
\begin{eqnarray*}
 S &=& {\Gamma(\sum\limits_{i=1}^\ell m_i)\over
  \prod\limits_{i=1}^\ell \Gamma(m_i)}  \sum\limits_{k=1}^{\infty}
   {\buildrel \begin{matrix}{0\le \sum\limits_{i=1}^{\ell-1}x_i\le1}\\
                           {0\le x_j \le 1}
             \end{matrix}
   \below {\int\int\cdots\int} }
   { (\prod\limits_{i=1}^{\ell-1} x_i^{m_i-1})
    (1-\sum\limits_{i=1}^{\ell-1} x_i)^{m_\ell-1} \over
    [k+a_\ell + \sum\limits_{i=1}^{\ell-1} x_i(a_i-a_\ell)]^{\sum\limits_{i=1}^\ell m_i}}
    \, dx_1 dx_2 \cdots dx_{\ell-1} \\ \  \\
 &=&  {1\over \prod\limits_{i=1}^\ell \Gamma(m_i)} \sum\limits_{k=1}^{\infty}
    {\buildrel \begin{matrix}{0\le \sum\limits_{i=1}^{\ell-1}x_i\le1}\\
                           {0\le x_j \le 1}
             \end{matrix}
   \below {\int\int\cdots\int} }
      dx_1 dx_2 \cdots dx_{\ell-1}\,
    (\prod\limits_{i=1}^{\ell-1} x_i^{m_i-1})
    (1-\sum\limits_{i=1}^{\ell-1} x_i)^{m_\ell-1} \times \\
    && \hspace{4cm} \times
    \int_0^{+\infty}dy \,  y^{(\sum\limits_{i=1}^\ell m_i-1)}
    e^{-y[k+a_\ell + \sum\limits_{i=1}^{\ell-1} x_i(a_i-a_\ell)]}
    \\ \  \\
    &=&  {1\over \prod\limits_{i=1}^\ell \Gamma(m_i)}
   {\buildrel \begin{matrix}{0\le \sum\limits_{i=1}^{\ell-1}x_i\le1}\\
                           {0\le x_j \le 1}
             \end{matrix}
   \below {\int\int\cdots\int} }
      dx_1 dx_2 \cdots dx_{\ell-1}\,
    (\prod\limits_{i=1}^{\ell-1} x_i^{m_i-1})
    (1-\sum\limits_{i=1}^{\ell-1} x_i)^{m_\ell-1} \times \\
    && \hspace{4cm} \times
    \int_0^{+\infty}  y^{(\sum\limits_{i=1}^\ell m_i-1)}
    ~ {e^{-y[a_\ell+1 + \sum\limits_{i=1}^{\ell-1} x_i(a_i-a_\ell)]}
    \over 1-e^{-y}} \, dy ~.
 \end{eqnarray*}
To pass from the first expression to the second, we used formula
(\ref{laplace2}). Wheelon does not present explicitly the above
formula; he only describes his method in a simple case and then he
lists exact results for series for which $P(k)$ is a product of up
to four distinct monomials. The result for Andreoli's problem
(which will be computed in the next section) does not appear in
\cite{Wheelon}.

\paragraph*{Andreoli's problem again}
When $m_i=1$,  for all $i$, Feynman's integral can be simplified
considerably. In particular, one can write
\begin{eqnarray*}
  {1\over \alpha_1 \alpha_2 \cdots \alpha_\ell } &=& (\ell-1)!
   \int_0^1 dx_1 \int_0^{x_1}dx_2 \cdots
   \int_0^{x_{\ell-2}}dx_{\ell-1}\times \\
  && {1 \over [\alpha_1+x_1(\alpha_2-\alpha_1)+\cdots+
                        x_{\ell-2}(\alpha_{\ell-1}-\alpha_{\ell-2})+
             x_{\ell-1}(\alpha_\ell-\alpha_{\ell-1})]^\ell}~,
\end{eqnarray*}
 If $\ell=n+1$ and $\alpha_i=k+i-1$, $i=1,2,\dots,n+1$, then
\begin{eqnarray}
  {1\over k (k+1)(k+2) \cdots (k+n) } &=& n!
  \int_0^1 dx_1 \int_0^{x_1}dx_2 \cdots
   \int_0^{x_{n-1}}dx_n
   {1\over
    [k+\sum\limits_{i=1}^n x_i]^{n+1}}  \nonumber \\
  &=& \int_0^1 dx_1 \int_0^1dx_2 \cdots
   \int_0^1 dx_n
   {1\over
    [k+\sum\limits_{i=1}^n x_i]^{n+1}} \nonumber \\
    &=& {1\over n!} \int_0^1 dx_1 \int_0^1dx_2 \cdots
   \int_0^1 dx_n
   \int_0^{+\infty}  y^n
    e^{-y[k+\sum\limits_{i=1}^n x_i]} \, dy \nonumber \\
  &=& {1\over n!}\int_0^{+\infty}  e^{-yk} y^n
      \lp\int_0^1  e^{-yx} \, dx \rp^n \, dy \nonumber \\
  &=& {1\over n!}\int_0^{+\infty}  e^{-yk}
      (1-e^{-y})^n \, dy ~,
      \label{eq:11}
\end{eqnarray}
and
\begin{eqnarray*}
  \sum_{k=1}^\infty {1\over k (k+1)(k+2) \cdots (k+n) }
  &=& {1\over n!}\int_0^{+\infty}
      (1-e^{-y})^n \sum_{k=1}^\infty e^{-yk} \, dy\\
  &=& {1\over n!}\int_0^{+\infty}
      (1-e^{-y})^{n-1} e^{-y} \, dy ~,
\end{eqnarray*}
which is identical to equation (\ref{eq:10}) found by the partial
fractions method.

\paragraph*{A Related Result}
 Equation (\ref{eq:11}) is valid for any positive number $k$ (we shall rename it $x$), not
 necessarily an integer. If instead we sum over $n$
 \begin{eqnarray}
 \sum_{n=0}^\infty {1\over x (x+1)(x+2) \cdots (x+n) }
  &=&  \int_0^1  u^{x-1} \,
      \sum_{n=0}^\infty {(1-u)^n\over n!} \, du
      \nonumber \\
  &=&  e\, \int_0^1  u^{x-1} e^{-u} \, du~,
 \label{eq:6}
 \end{eqnarray}
 where we have made a change of variable $u=e^{-x}$. In this expression
 we can now expand $e^{-u}$ in a Taylor series and perform
 the elementary integrals to find
 \begin{equation}
 \sum_{n=0}^\infty {1\over x (x+1)(x+2) \cdots (x+n) }
 \= e \, \lb {1\over x} -{1\over 1!}\,{1\over x+1}
                        +{1\over2!}\,{1\over x+2}-\cdots\rb~.
 \label{eq:5}
 \end{equation}

\paragraph*{Additional Results}

The series proposed by Andreoli and its generalization
(\ref{eq:9}) are well-known: see page 234 of \cite{GR}. Equations
(\ref{eq:5}) and (\ref{eq:6}) are also found in \cite{GR} (page
317). Equation (\ref{eq:5})---although not equation
(\ref{eq:6})---is presented as an exercise (ex. 126, page 272) in
Knopp's book \cite{Knopp}. In this section we extend the results
presented above to other series. The results obtained seem to be
new; we have not found them in  \cite{GR} or \cite{Knopp}.

Using the Feynman integral and equation (\ref{laplace2}), we can
show that
\begin{eqnarray}
  {1\over k (k+2)(k+4) \cdots (k+2n) }
  &=& {1\over n!}\int_0^1 u^{k-1} \, (1-u^2)^{n-1}  \, du ~.
\label{eq:new1}
\end{eqnarray}
By renaming $k$ as $x$, summing over $n$, and following the steps
that led to equation (\ref{eq:5}), we can find
\begin{eqnarray}
  \sum_{n=0}^\infty {1\over x (x+2)(x+4) \cdots (x+2n) }
  &=& e^{1/2} \,\int_0^1 u^{x-1} \, e^{-u^2/2}  \, du
  \nonumber \\
  && \hspace{-2cm}
  = e^{1/2} \, \lb {1\over x} -{1\over 2\,1!}\,{1\over x+2}
                        +{1\over 2^2 2!}\,{1\over x+4}-\cdots\rb~.
\label{eq:new2}
\end{eqnarray}
On the other had, summing in (\ref{eq:new1}) over $k$ gives:
\begin{eqnarray*}
  \sum_{k=1}^\infty {1\over k (k+2)(k+4) \cdots (k+2n) }
  &=& {1\over 2^nn!}\int_0^1  \, (1-u^2)^{n} \,{1\over 1-u} \, du ~.
\end{eqnarray*}
We could have arrived at the same integral representation if we
had used the partial fraction decomposition. The integral of the
right hand side may be rewritten as
$$
  \int_0^1  \, (1-u^2)^{n-1} \,(1+u) \, du
  ~=~
  \int_0^1  \, (1-u^2)^{n-1} \,du +
  \int_0^1  \, (1-u^2)^{n-1} \, u\, du~.
$$
These integrals in turn can be computed\footnote{Actually, the
second integral is  basic.} by the use of the B-function
(\ref{eq:beta}).
Finally, we find
\begin{eqnarray}
  \sum_{k=1}^\infty {1\over k (k+2)(k+4) \cdots (k+2n) }
  &=& {1\over 2n} \lp {1\over(2n-1)!!}+{1\over(2n)!!}\rp ~.
\label{eq:new4}
\end{eqnarray}

Empowered with the results that have been presented so far, one
can ask if we can prove similar results for the series
\begin{eqnarray*}
 F(x) &=& \sum_{n=0}^\infty {1\over x (x+\ell)(x+2\ell) \cdots (x+n\ell)
 }~, \\
 S(\ell) &=& \sum_{k=1}^\infty {1\over k (k+\ell)(k+2\ell) \cdots (k+n\ell)
 }~.
\end{eqnarray*}
Indeed, it is straightforward to show that
\begin{eqnarray}
 F(x) &=& e^{1/\ell} \,\int_0^1 u^{x-1} \, e^{-u^\ell/\ell}  \, du
 \nonumber \\
  &=& e^{1/\ell} \, \lb {1\over x} -{1\over \ell\,1!}\,{1\over x+\ell}
   +{1\over \ell^2 2!}\,{1\over x+2\ell}-\cdots\rb~,
\label{eq:new3}
\end{eqnarray}
and that
\begin{eqnarray*}
 S(\ell) &=& {1\over \ell^n n!}\, \int_0^1
 (1-u^\ell)^n\,{1\over 1-u}\, du ~.
\end{eqnarray*}
The integral of the right hand side may be rewritten as
$$
  \int_0^1  \, (1-u^\ell)^{n-1} \,(1+u+u^2+\dots+u^{\ell-1}) \, du
  ~,
$$
that is, a sum of integrals of the form (\ref{eq:beta}). We thus
find
\begin{eqnarray}
 S(\ell) &=& {1\over \ell\, n}\, \lp
 {1\over(\ell n)!^\ell}+{1\over(\ell n-1)!^\ell}+
 {1\over(\ell n-2)!^\ell}+\dots+{1\over(\ell n-(\ell-1))!^\ell}\rp
 ~,
 \label{eq:new6}
\end{eqnarray}
where the symbol
 $!^\ell$ 
 is an extension of the double
factorial and stands for the product of all numbers between 1 and
$\ell n-m$ that, when divided by $\ell$, give residue $\ell-m$.
For example, $(3n)!!!=(3n)!^3=3\cdot6\cdot9\cdots(3n)$,
$(3n-1)!!!=(3n-1)!^3=2\cdot5\cdot8\cdots(3n-1)$, and so on.

\paragraph*{Conclusion}

Taking a problem proposed by Andreoli as the origin of this
article, we have attempted to draw the attention of the reader to
the use of  integral transforms in the computation of exact values
for series and make him  appreciate the ease and transparency of
the method. In the process, we have derived some results
(\ref{eq:new3}) and (\ref{eq:new6}) that seem to be new or, at
least, do not appear in the commonly used references \cite{GR} and
\cite{Knopp}.

The method can be applied to a wide variety of series as explained
in the articles of Efthimiou \cite{Efthimiou,Efthimiou2} and Lesko
and Smith \cite{Lesko} and the book of Wheelon \cite{Wheelon}. The
reader is encouraged to look at these references for additional
information. Hopefully,  he will use the method to search for new
results. He may begin with the series:
\begin{eqnarray*}
 F(x;a,b) &=& \sum_{n=0}^\infty {1\over [a+xb] [a+(x+\ell)b]
                  [a+(x+2\ell)b] \cdots [a+(x+n\ell)b] }~, \\
 S(\ell;a,b) &=& \sum_{k=1}^\infty {1\over [a+kb] [a+(k+\ell)b]
                  [a+(k+2\ell)b] \cdots [a+(k+n\ell)b] }~.
\end{eqnarray*}

\paragraph*{Acknowledgements}

Few days before the posting of the present paper on the archives,
I received a message from M. L. Glasser drawing my attention to
one of his papers \cite{Glasser} in which he had extended
Wheelon's method as described in his original paper
\cite{Wheelon2} to additional series. Glasser's paper has
considerable overlapping with \cite{Efthimiou}, \cite{Lesko} and
\cite{Efthimiou2} although the presentation differs from these
papers.


\end{document}